\documentclass[10pt]{article}
\usepackage{amssymb}
\usepackage{amsmath}
\usepackage{amsthm}
\theoremstyle{plain}
\newtheorem{theorem}{Theorem} \newtheorem{lemma}{Lemma}[section]
\newtheorem{propo}{Proposition}[section]

\newtheorem{corol}{Corollary}[section] \newtheorem{defin}{Definition}[section]

  \newcommand{\ep}{\varepsilon}
  
 \newcommand{\e}{\ep} 
\newcommand{\N}{\mathbb{N}} 
\newcommand{\R}{\mathbb{R}} \newcommand{\Z}{\mathbb{Z}}
 \newcommand{\C}{\mathbb{C}}
\newcommand{\F}{\mathbb{F}}

\newcommand{\Mat}{\mbox{Mat}\,} 
\newcommand{\Id}{\mbox{Id}}\newcommand{\Ima}{\mbox{Im}\,}
\newcommand{\rk}{\mbox{rk}} 
\newcommand{\rank}{\mbox{Rank}\,} \newcommand{\dck} {\dim_{\cK}\,}
\newcommand{\Span}{\mbox{Span}\,}\newcommand{\End}{\mbox{End}\,}
\newcommand{\Ker}{\mbox{Ker}\,}
\newcommand{\Tr}{\mbox{Tr}}\newcommand{\Fix}{\mbox{Fix}}
\newcommand{\cM}{\mathcal{M}}\newcommand{\cD}{\mathcal{D}}
\newcommand{\cN}{\mathcal{N}} 
\newcommand{\cH}{\mathcal{H}}
\newcommand{\cK}{\mathcal{K}}
\newcommand{\cF}{\cK\F_r}
  
\newcommand{\mnk}{\Mat_{n_k\times n_k}(\cK)}
\newcommand{\mnd}{\Mat_{n\times n}(D)}
\newcommand{\hp} {\hat{\phi}} \newcommand{\hh}{\hat{\theta}}
 
\newcommand{\dom} {\mbox{dom\,}}
\newcommand{\ul} {\underline{\lambda}}
\title{Convergence and limits of linear representations of
finite groups
\footnote{
\, Research partly sponsored by MTA Renyi ``Lendulet'' Groups
and Graphs Research Group}}

\author{G\'abor Elek}

\begin{document}
\maketitle

\begin{abstract}
Motivated by the theory of graph limits, we introduce and 
study the convergence and limits of linear representations
of finite groups over finite fields. The limit objects are
infinite dimensional representations of free groups in continuous
algebras. We show that under a certain integrality condition, the
algebras above are skew fields. Our main result is the extension of
Schramm's characterization
of hyperfiniteness to linear representations. 
\end{abstract}
\footnotesize
\tableofcontents
\normalsize
\newpage
\section{Introduction}
\fontsize{11pt}{13pt}\selectfont
\vskip 0.2in
\noindent{\bf Unitary representations.}
Our primary motivation and model example is the view of infinite dimensional
unitary representations into tracial von Neumann algebras as limits of
finite dimensional unitary representations.
By a finite dimensional unitary representation (of degree $r$),
we mean a homomorphism
$\kappa:\F_r\to U(n)$ of the free group on $r$ generators 
into the unitary group
$U(n)$. Note that such representations can be given by the $r$-tuple
$\{\kappa(\gamma_i\}^r_{i=1}$, where $\{\gamma_i\}^r_{i=1}$ are the standard 
generators of the free group $\F_r$. We say that a sequence $\{\kappa_k:
\F_r\to U(n_k)\}^\infty_{k=1}$
of finite dimensional unitary represenations are {\it convergent} if
for any $w\in\F_r$
$$\Tr(w)=\lim_{k\to\infty}\Tr_{n_k} (\kappa_k(w))$$
exists, where $\Tr_{n_k}$ stands for the normalized trace function on
the complex matrix algebra $\Mat_{n_k\times n_k}(\C)$.
Note that from each sequence of representations one can choose a convergent
subsequence. The limit of a convergent system $\{\kappa_k:
\F_r\to U(n_k)\}^\infty_{k=1}$
is defined as a representation of $\F_r$ into the unitary group of a tracial
von Neumann algebra via the $GNS$-construction, given below.

The function $\Tr$ extends to the group algebra $\C\F_r$ as a trace, that is
a linear functional satisfying
$$\Tr(ab)=\Tr(ba)\,.$$
Let $I\subset \C\F_r$ be the set of elements $a$, for which $\Tr(a^*a)=0$.
It is not hard to see that $I$ is an ideal of $\C\F_r$ and the trace function
$\Tr$ vanishes on $I$. Let $A=\C\F_r/I$, then
$$<[p],[q]>=\Tr[pq^*]$$
is well defined and gives rise to an inner product structure on the algebra
$A$. Let $\cH$ be the Hilbert completion of $A$. Then the
left multiplication $L_{[p]}:A\to A$ defines
a bounded linear representation of $A$ on $\cH$. The weak closure of the
image of $A$ is a tracial von Neumann algebra $\cN$ equipped with a trace (the
extension of $\Tr$). Also, we have a natural homomorphism $S:\F_r\to U(\cN)$
such that for any $w\in \F_r$, 
$$\Tr(S(w))=\lim_{k\to\infty} \Tr_{n_k}(\kappa_{k}(w)).$$
Thus $S$ can be viewed as the limit object of the finite dimensional unitary
representations $\{\kappa_{k}\}^\infty_{k=1}$. One can ask the following question. If
$S:\F_r\to U(\cN)$ is a representation of the free group into the unitaries
of a tracial von Neumann algebra, is it true that $S$ is the 
limit of finite dimensional
unitary representations. This question is called the Connes Embedding Problem
\cite{Pestov}.
\vskip 0.2in
\noindent
{\bf Graph limits.} We are also strongly motivated by the emerging theory
of graph limits. The original definition of graph convergence is due to
Benjamini and Schramm \cite{Benjamini} 
(see also the monography of Lov\'asz \cite{Lovasz}).
Here we consider the limit of Schreier graphs (see e.g. \cite{Elek2012}).
Let $\lambda:\F_r\to S(n)$ be a homomorphism of the free group into the 
symmetric group of permutations on $n$ elements. These homomorphism
are in bijective correspondence with Schreier graphs. The vertex set of the
corresponding
graph $G_\lambda$ is the set $[n]=\{1,2,\dots,n\}$. Two vertices $a$ and $b$ are
connected with an edge labeled by the generator $\gamma_i$, if
$\lambda(\gamma_i)(a)=b$. A sequence of permutation representations (or Schreier graphs)
$\{\lambda_k:\F_r\to S(n_k)\}^\infty_{k=1}$
is called {\it convergent} if for any $m\geq 1$ and 
$m$-tuple
$\{w_1,w_2,\dots,w_m\}\subset \F_r$

$$\lim_{k\to\infty} \frac{ |\Fix_k(\lambda_k(w_1))\cap \Fix_k(\lambda_k(w_2))\cap
\dots\cap \Fix_k(\lambda_k(w_m))|} {n_k}$$
exists, where $\Fix_k(\lambda_k(w))$ is the fixed point set of the permutation
$\lambda_k(w)$. Note that the original definition is somewhat different from
the one above, nevertheless a simple inclusion-exclusion argument 
shows that the
two definitions are equivalent. One can define various limit objects for
such convergent sequences e.g. the invariant random subgroups 
(see \cite{Abert}).
A notion of limit, analogous to the unitary case, can be defined the
following way. Let $(X,\mu)$ be a standard Borel
probability measure space equipped with a countable measure preserving
 equivalence relation $E$
(see \cite{Kechris}).
An $E$-transformation is a measure preserving bijection $T:X\to X$ such
that for almost all $p\in X$, $p\equiv_E T(p)$. A zero transformation
is an $E$-transformation $S$ such that $\mu(\Fix(S))=1$.
The full group of $E$, $[E]$ is the quotient of the
group of $E$-transformations by the 
normal subgroup of zero transformations. Note that if $Q\in [E]$, then
$\Fix(Q)$ is well-defined up to a zero measure perturbation. We call
a homomorphism $\lambda:\F_r\to [E]$ {\it generating} if for almost
all $p\in X$: for any $q$ such that $p\equiv_E q$, there
exists $w\in \F_r$ such that
$\lambda(w)(p)=q$.
A generating representation $\lambda:\F_r\to [E]$ is a limit of the convergence
system of permutation representations $\{\lambda_k\}^\infty_{k=1}$
if for any $m\geq 1$ and 
$m$-tuple
$\{w_1,w_2,\dots,w_m\}\subset \F_r$

$$\lim_{k\to\infty} \frac{ |\Fix_k(\lambda_k(w_1))\cap 
\Fix_k(\lambda_k(w_2))\cap
\dots\cap \Fix_k(\lambda_k(w_m))|} {n_k}= $$ $$=
\mu\left(\Fix(\lambda(w_1)) \cap \Fix(\lambda(w_2)) \cap
\dots \cap \Fix(\lambda(w_m))\right)\,. $$

\noindent
Again, for each convergent sequence $\{\lambda_k\}^\infty_{k=1}$ one can
find a limit representation into some full group. On the other hand,
it is not known, whether any generating representation $\lambda:\F_r\to [E]$
can be obtained as a limit. 

\vskip 0.2in
\noindent
{\bf Schramm's Theorem.} Let $\lambda:\F_r\to S(n)$ be a permutation
representation with corresponding Schreier graph $G_\lambda([n],E_\lambda)$.
We say that a convergence sequence of representations 
$\{\lambda_k:\F_r\to S(n_k)\}^\infty_{k=1}$
is hyperfinite if for any $\epsilon>0$, there exists an integer
$K_\epsilon>0$ such that for any $k\geq 1$, one can remove $\e n_k$ edges 
from the graph $G_{\lambda_k}$ in such a way, that all the components of the
remaining graph have at most $K_\epsilon$ vertices. Schramm \cite{Schramm}
proved that the hyperfiniteness of the sequence is equivalent to the
amenability of its limit (see also \cite{Elek2012}). 
A generating representation $\lambda:\F_r\to [E]$
is amenable if $E$ is a hyperfinite (amenable) equivalence relation 
\cite{Kechris}.

\section{Definitions and results}
In the course of our paper we fix a finite field $\cK$. Our goal is to study
the convergence and limits of finite dimensional representations $\theta:
\F_r\to \Mat_{n\times n}(\cK)$. Note that such representations are in one to one
correspondence with injective
linear representations $\pi:\Gamma\to GL(n,\cK)$, where
$\Gamma$ is a finite group of $r$ marked generators.
\begin{defin} \label{crucial}
A sequence of finite dimensional representations \\
$\{\theta_k:\F_r\to \Mat_{n_k\times n_k}(\cK)\}^\infty_{k=1}$
is convergent if for all $n\geq 1$ and matrix $A\in \Mat_{n\times n}(\cK\F_r)\,,$
$$\lim_{k\to\infty} \rk_{n_k}^n(\theta_k(A))$$
exists, where $\theta_k(A)$ is the image of $A$ in 
$$\Mat_{n_k n\times n_k n}(\cK)\cong\Mat_{n\times n}(Mat_{n_k\times n_k}(\cK))$$ and 
$$\rk_{n_k}^n(\theta_k(A))=\frac{\rank (\theta_k(A))}{n_k}\,.$$
\end{defin}

\noindent
 Note that $\theta_k$
naturally extends to the group algebra $\cF$ and we denote the extension
by $\theta_k$, as well.
We will make clear at the end of
 Section \ref{convergent}, why we consider matrices instead of
single elements of the group algebra $\cF$. 
Now we define the limit objects for convergent
sequences. The objects we need, continuous algebras, were defined by John von 
Neumann in the thirties \cite{Neumann}. 
Let $R$ be a separable, continuous $\cK$-algebra with rank function
$\rk_R$ (see Section \ref{cont} for a brief survey on continuous algebras).
\begin{defin}
Let $\{\theta_k\}^\infty_{k=1}$ be finite dimensional representations as above.
A representation $\theta:\F_r\to R$, that is a homomorphism of the free group
into the group of invertible elements of the continuous algebra $R$ is
a limit of $\{\theta_k\}^\infty_{k=1}$, if for any $n\geq 1$ and
$A\in \Mat_{n\times n}(\cF)\,$
$$\lim_{k\to\infty} \rk^n_{n_k}(\theta_k(A))=\rk^n_R(\theta(A))\,,$$
where $\rk^n_R$ is the matrix rank on $\Mat_{n\times n}(R)$. 
\end{defin}
\noindent
Note, that if $\cN$ is a tracial von Neumann algebra, then $\cN$ is
equipped with a natural rank function and its completion is a continuous
rank regular ring; the algebra of affiliated operators \cite{Thom}.
Hence the limit of finite
dimensional unitary representations can also be viewed as a homomorphism
into the group of invertible elements of a continuous algebra.
\noindent
Our first theorem is about the existence of limits.
\begin{theorem} \label{elso}
For any convergent sequence of finite dimensional representations
$\{\theta_k\}^\infty_{k=1}$, there exists a separable, continuous
$\cK$-algebra $R$ and a representation $\theta:\F_r\to R$ such that
$\theta$ is the limit of $\{\theta_k\}^\infty_{k=1}$.
\end{theorem}
\noindent
It turns out that under a certain integrality condition the limit is
unique. We say that the convergence sequence $\{\theta_k\}^\infty_{k=1}$
satisfies the Atiyah condition, if for any $n\geq 1$ and
 $A\in \Mat_{n\times n}(\cK\F_r)$
$$\lim_{k\to\infty} \rk^n_{n_k}(\theta_k(A))\in \Z\,.$$
The condition above is intimately related to Atiyah's Conjecture
 on the integrality of
the $L^2$-betti numbers (see \cite{Linnell} and \cite{Linnell2}).
\begin{theorem} \label{masodik}
If the convergence sequence of linear representations 
$\{\theta_k\}^\infty_{k=1}$ satisfies
the Atiyah condition,
then there exists a skew field $D$ over the base field $\cK$ and
a homomorphism $\theta:\F_r\to D$ (that is a
homomorphism into the multiplicative group of non-zero elements of $D$) 
such that $\theta$ is the limit
of $\{\theta_k\}^\infty_{k=1}$ and $\Ima (\theta)$ generates $D$. Moreover,
if $\theta':\F_r\to D'$ is another generating 
limit homomorphism into a skew field
$D'$, then there exists a skew field isomorphism $\pi:D\to D'$ such that
$\pi\circ \theta=\theta'$.
\end{theorem}
\noindent
If the Atiyah condition is satisfied, we will be able to 
generalize Schramm's Theorem cited in
the Introduction. It is worth to note that hyperfinite sequences of 
graphs are basically the opposites
of expander sequences. The notion of expander sequences 
for linear representations
were introduced by Lubotzky and Zelmanov (see also \cite{Dvir}). We say that
a sequence of linear representations 
$\{\theta_k:\F_r\to \Mat_{n_k\times n_k}(\cK)\}^\infty_{k=1}$ is a 
{\it dimension expander} if there exists $\alpha>0$ such
 that for all $k\geq 1$ and linear
subspace $W\subset \cK^{n_k}$ with $\dim_{\cK}(W)\leq \frac{n_k}{2}$
$$\dim_{\cK} (W+\sum^r_{i=1}\theta_k(\gamma_i)(W))\geq (1+\alpha) \dck(W)\,,$$
where $\{\gamma_i\}^r_{i=1}$ are the standard generators of $\F_r$.
Note that 
$$\sup_{W,\, \dck(W)\leq \frac{n_k}{2}} \frac{\dim_{\cK} 
(W+\sum^r_{i=1}\theta_k(\gamma_i)(W))}
{\dck(W)}$$
is the linear analogue of the Cheeger constant of a graph. 
It was observed in \cite{Dvir} that
a random choice of the $r$-tuple $\{\theta_k(\gamma_1),
\theta_k(\gamma_2),\dots,\theta_k(\gamma_r)\}^\infty_{k=1}$ leads
 to a dimension expander with
probability one, provided that $r$ is large enough. Later, 
Bourgain and Yehudayoff \cite{Bourgain} constructed
explicit families of dimension expanders. Using the 
linear graph theory vocabulary:
{\it subsets}$\to${\it linear subspaces}, 
{\it disjoint}$\to$ {\it independent}, {\it cardinality}$\to$
{\it dimension}, we can define the hyperfiniteness 
for sequences of linear representations.
\begin{defin}
The linear representations \\
$\{\theta_k:\F_r\to \Mat_{n_k\times n_k}(\cK)\}^\infty_{k=1}$ form
a hyperfinite sequence if for any $\epsilon>0$ there 
exists $K_\epsilon>0$ such that
for all $k\geq 1$, we have $\cK$-linear 
subspaces $V^k_1, V^k_2,\dots, V^k_{t_k}\subset \cK^{n_k}$ such that
\begin{itemize}
\item For any $1\leq j\leq t_k$, $\dck(V^k_j)\leq K_\epsilon$.
\item $\{V^k_j+\sum^r_{i=1}\theta_k(\gamma_i)V^k_j\}^{t_k}_{j=1}$
 are independent subspaces such that
$$\dim_{\cK} (V^k_j+\sum^r_{i=1}\theta_k(\gamma_i)V^k_j)<
(1+\epsilon) \dim_{\cK} (V^k_j).$$
\item $\sum^{t_k}_{j=1}\dck(V^k_j)\geq (1-\epsilon) n_k\,.$
\end{itemize}
\end{defin}
\noindent
Our main result is the generalization of Schramm's Theorem for
convergent sequences of linear representations. 
\begin{theorem} \label{harmadik}
Let $\{\theta_k\}^\infty_{k=1}$ be a convergent sequence of 
linear representations satisfying the
Atiyah condition. Let $\theta:\F_r\to D$ be the unique limit representation of
$\{\theta_k\}^\infty_{k=1}$.
 Then  $\{\theta_k\}^\infty_{k=1}$
is hyperfinite if and only if $D$ is an amenable skew field.
\end{theorem}

\section{Universal localizations} \label{local}
In this section we recall the notion of
universal localization from the book of Cohn \cite{Cohn}.
Let $R$ be a unital ring, $\Sigma$ be a set of matrices
over $R$ and $f:R\to S$ be a unital ring homomorphism.
Let $\Sigma^f$ be the image of $\Sigma$ under $f$.
If $\Sigma$ is the set of matrices whose images under $f$ are invertible
, then $R^f(S)$ denotes the set of entries
in the inverses $M^{-1}$, for $M\in\Sigma^f$. 
We call $R^f(S)$ the rational closure of $R$ in $S$. According to
Theorem 7.1.2 \cite{Cohn} $R^f(S)$ is a ring containing $\Ima(f)$.
An important tool for the understanding of the ring $R^f(S)$ is the following
 variant
of Cramer's Rule.
\begin{lemma}[Proposition 7.1.3 \cite{Cohn}] \label{cramer}
For any element $p\in R^f(S)$, there exists $n\geq 1$
and $Q\in Mat_{n\times n}(\Ima(f))$, $A\in \Mat_{n\times n}(\Ima(f))$ invertible
in $\Mat_{n\times n}(S)$ and $B\in \Mat_{n\times n}(R^f(S))$ in the form of
$B=\begin{pmatrix}I & u \\ 0 & 1 \end{pmatrix}$ such that
$$Q=A \begin{pmatrix}I & 0 \\ 0 & p \end{pmatrix} B\,.$$
\noindent
Note that $I$ denotes the unit matrix of size $n-1$.
\end{lemma}
\noindent
Recall that the division closure of $\Ima(f)$ in $S$ is the smallest
ring in $S$ containing $\Ima(f)$ closed under taking inverses.
The following lemma is given as an exercise in \cite{Cohn}.
\begin{lemma} \label{division}
The division closure $D(f)$ is a subring of $R^f(S)$.
\end{lemma}
\proof
It is enough to prove that if $p\in R^f(S)$ is invertible in $S$, then
$p^{-1}\in R^f(S)$. By Lemma \ref{cramer}, we can write
$$Q=A \begin{pmatrix}I & 0 \\ 0 & p \end{pmatrix} B\,.$$
\noindent
So
$\begin{pmatrix} I & 0 \\ 0 & p \end{pmatrix}=A^{-1}Q B^{-1}$
that is $\begin{pmatrix} I & 0 \\ 0 & p^{-1} \end{pmatrix}=B Q^{-1} A\,.$
Since all the entries of $B,Q^{-1}$ and $A$ are, by definition, in the
subring $R^f(S)$, we have that $p^{-1}\in R^f(S)$. \qed 

\vskip 0.2in
\noindent
Now, let $\Sigma$ be a set of square matrices over $R$. The universal
localization of $R$ with respect to $\Sigma$ is the unique ring $R_{\Sigma}$ 
equipped with a homomorphism
$\lambda:R\to R_\Sigma$ such that the elements
of $\Sigma^\lambda$ are all invertible matrices and if $f:R\to S$ is
an arbitrary homomorphism and the elements of $\Sigma^f$ are all invertible
matrices, then
there exists a unique homomorphism $\overline{f}:R_\Sigma\to S$ such that
$\overline{f}\circ \lambda=f$. Let $D$ be a skew field and $f:R\to D$
be a homomorphism. We call $D$ epic if $\Ima(f)$ generates $D$ as a skew field. 
\begin{propo} [Theorem 7.2.2 \cite{Cohn}]
If $D$ is epic and $\Sigma_f$ is the
set of matrices over $R$ whose images in $D$ are invertible, then the 
universal localization
$R_{\Sigma_f}$ is a local ring with residue-class field isomorphic to $D$.
\end{propo}
\noindent
We have the following corollary.
\begin{corol} \label{korolla}
If $f_1:R\to D_1$ and $f_2:R\to D_2$ are two epic maps and
for all $n\geq 1$ and for all matrices $A\in \Mat_{n\times n}(R)$
$$\rk_{D_1}(f_1(A))=\rk_{D_2}(f_2(A)),$$
then $D_1 \cong D_2$ and there is an isomorphism $\iota:D_1\to D_2$ such 
that 
\begin{equation} \label{e524} \iota\circ f_1=f_2 .\end{equation}
\end{corol}
\proof
By our condition, $\Sigma_{f_1}=\Sigma_{f_2}=\Sigma$. Let
$\lambda':R\to R_{\Sigma}/M$ be the natural map, where $M$ is the unique maximal
ideal of $R_{\Sigma}$. Then there exist two skew field isomorphisms
$\pi_1:R_{\Sigma}/M \to D_1$ and $\pi_2:R_{\Sigma}/M\to D_2$ such that
$\pi_1\circ \lambda'=f_1$ and $\pi_2\circ \lambda'=f_2$. Then we can choose
$\iota=\pi_2\circ\pi^{-1}$ to satisfy (\ref{e524}). \qed

\section{Continuous algebras} \label{cont}
In this section we recall the notion of a continuous algebra from the
book of Goodearl \cite{Goodearl} and present some important
examples.
A ring $R$ is called von Neumann regular if for any $a\in R$ there
exists $x\in R$ such that $axa=a$. In other words, $R$ is von Neumann regular
if any finitely generated left ideal is generated by a single idempotent.
A rank regular ring is a unital regular ring $R$ equipped with a rank function
$\rk_R$ satisfying the conditions below.
\begin{itemize}
\item $0\leq \rk_R(a)\leq 1\,,$ for any $a\in R$.
\item $\rk_R(a)=0$ if and only if $a=0$.
\item $\rk_R(1)=1.$
\item $\rk_R(a+b)\leq \rk_R(a)+\rk_R(b)\,.$
\item $\rk_R(ab)\leq \rk_R(a),\rk_R(b)$
\item $\rk_R(e+f)=\rk_R(e)+\rk_R(f)$ if $e$ and $f$ are orthogonal idempotents.
\end{itemize}
Note, that in a rank regular ring an element is invertible if and only if
it has rank one. Also, a rank regular ring is a metric space with respect
to the distance function
$$d_R(a,b):=rk_R(a-b)\,.$$
\noindent
If a rank regular ring $R$ is complete with respect to the distance function,
then $R$ is called a {\it continuous ring}. We are particularly interested
in continuous $\cK$-algebras. The simplest examples are skew fields over $\cK$
and 
matrix rings over such skew fields. For these continuous algebras the rank
function may take only finitely many values. Another important example is due
to John von Neumann.
Let us consider the diagonal maps
$$d_n:\Mat_{2^n\times 2^n}(\cK)\to \Mat_{2^{n+1}\times 2^{n+1}}(\cK)\,.$$
\noindent
The maps preserve the normalized rank functions, hence the direct
limit $\varinjlim \Mat_{2^n\times 2^n}(\cK)$ is a rank regular ring. Its metric
completion $A_{\cK}$ is a continuous $\cK$-algebra, with a rank function $\rk_A$
taking all values in between zero and one.

\noindent
Note that if $R$ is a rank regular ring, the metric completion of $R$
is always a continuous ring \cite{Halperin1}. Also, if $R$ is a rank regular
ring, then for any $n\geq 1$ the matrix ring $\Mat_{n\times n}(R)$ can be
 equipped
with a unique matrix rank function $\rk^n_R$ such that $\rk^n_R(\Id)=n$ and
\begin{equation}
\rk^n_R \begin{pmatrix} A & 0 \\ 0 & B \end{pmatrix}=
\rk^k_R(A) +\rk^{n-k}_R(B) \end{equation}
if $A\in \Mat_{k\times k}(R), B\in \Mat_{(n-k)\times (n-k)}(R)\,$ \cite{Halperin2}.
Finally, let us recall the notion of the ultraproduct of finite dimensional
matrix algebras. This construction will be crucial in our paper.
Let $M=\{\Mat_{n_k\times n_k}(\cK)\}^\infty_{k=1}$ be a sequence of matrix algebras
over our base field $\cK$ equipped with the normalized rank
functions $\{\rk_{n_k}\}^\infty_{k=1}$ such that $n_k\to \infty$. 
Let $\omega$ be an
ultrafilter on the natural numbers
and let $\lim_\omega$ be the associated ultralimit. 
The ultraproduct $\cM_M$ of the algebras 
$\{\Mat_{n_k\times n_k}(\cK)\}^\infty_{k=1}$
can be defined the following way.
Consider the elements
$$\{a_k\}^\infty_{k=1}\in \prod^\infty_{k=1}\Mat_{n_k\times n_k}(\cK)\,,$$
for which $\lim_{\omega} \rk_{n_k}(a_k)=0\,.$
It is easy to see that these elements form an ideal $I_M$. The ultraproduct
is defined as
$$\cM_M:=\prod^\infty_{k=1}\Mat_{n_k\times n_k}(\cK)/I_M\,.$$
The $\cK$-algebra $\cM_M$ is a simple
continuous algebra equipped with a rank function \cite{ElekSzabo}
$$\rk_\cM \left([\{a_k\}^\infty_{k=1}]) \right)=\lim_\omega \rk_{n_k}(a_k)\,,$$
where $[\{a_k\}^\infty_{k=1}]$ denotes the class of $\{a_k\}^\infty_{k=1}
\in \prod^\infty_{k=1}\Mat_{n_k\times n_k}(\cK)$
in $\cM_M$.
\section{Sofic algebras}
Let $R$ be a countable $\cK$-algebra over our finite base field $\cK$,
with $\cK$-linear basis $\{1=r_1,r_2,r_3,\dots\}$. Following 
Arzhantseva and Paunescu \cite{Arzhantseva}, we
call $R$ (linearly) sofic if there exists a function
$j:R\backslash \{0\}\to \R^+$ and a sequence of positive numbers $s_i\to 0$ 
such that for any $i\geq 1$ there exists $n_i\geq 1$ and a $\cK$-linear
map $\phi_i:R\to \Mat_{n_i\times n_i}(\cK)$ satisfying the conditions below:
\begin{itemize}
\item $\phi_i(1)=\Id$
\item $\rk_{n_i}(\phi_i(a))\geq j(a)$ if $0\neq a\in \Span\{r_1,r_2,\dots,r_i\}$
\item $\rk_{n_i}(\phi_i(ab)-\phi_i(a)\phi_i(b))<s_i$ if
$a,b\in \Span\{r_1,r_2,\dots,r_i\}$.
\end{itemize}
Such a system is called a {\it sofic representation} of $R$. Clearly, the
soficity of an algebra does not depend on the particular choice of the
basis $\{r_1,r_2,\dots\}$.
\noindent
Using the maps above, we can define a map $\phi:R\to \cM_M$, by
$$\phi(s):=[\{ \phi_i(s)\}^\infty_{i=1}]\,,$$
where $M=\{\Mat_{n_i\times n_i}(\cK)\}^\infty_{i=1}\,.$
By our assumptions, $\phi$ is a unital embedding.
Conversely, we have the following proposition. Note that the proposition was
already implicite in \cite{Arzhantseva}, the proof below was suggested
by the referee. 
\begin{propo} \label{soficity}
Let $R$ be a countable algebra over our base field $\cK$.
If $R$ can be embedded into an ultraproduct $\cM_M$, then $R$ is sofic.
\end{propo}
\proof Let $\{1=r_1,r_2,\dots\}$ be a basis for $R$ and let \\ $j(a):=
\frac{1}{2} \rk_{\cM}(\phi(a))\,.$
It is enough to prove that for any $\epsilon>0$ and $i\geq 1$,
there exists an integer $n\geq 1$ and a linear, unit preserving map
$\sigma:R\to \Mat_{n\times n}(\cK)$ such that
\begin{itemize}\item
$\rk_n(\sigma(a))\geq j(a)$ if $0\neq a\in\Span\{r_1,r_2,\dots, r_i\}\,.$
\item
$\rk_n(\sigma(ab)-\sigma(a)\sigma(b))<\e$ if 
$a,b\in\Span\{r_1,r_2,\dots, r_i\}\,.$
\end{itemize}
\noindent
Let $\phi:R\to\cM_M$ be the embedding. Lift $\phi$ to a unital linear map
$$\prod^\infty_{k=1}\phi_k: R\to \prod_k^\infty \Mat_{n_k\times n_k}(\cK)$$
\noindent 
by lifting the basis of $R$ first and then extending to $R$ linearly.
Then, for any $a,b\in R$ and $\e>0$, the set
$$\{k\in\N\,\mid\, \rk_{n_k}(\phi_k(ab)-\phi_k(a)\phi_k(b))<\e\}\in\omega\,.$$
\noindent
Also, for any $0\neq a\in R$ the set
$$\{k\in \N\,\mid\,\rk_{n_k}(\phi_k(a))>
\rk_{\cM}(\phi(a))/2\}\in\omega\,.$$
Hence, for any $i$ and $\e>0$, the set of all $k\in\N$ satisfying
$\rk_{n_k}(\phi_k(ab)-\phi_k(a)\phi_k(b))<\e$ for all 
$a,b\in\Span\{r_1,r_2,\dots, r_i\}\,$ and
$\rk_{n_k}(\phi_k(a))>\rk_{\cM}(\phi(a))/2$ for all $0\neq a\in 
\Span\{r_1,r_2,\dots, r_i\}\,$ is in the ultrafilter $\omega$, in
particular, it is not empty. \qed

\section{Amenable skew fields are sofic}
In this section we recall the notion of amenability for skew fields from 
\cite{Elek2006}. Let $D$ be a countable dimensional skew field over a finite
 base field $K$. We say that $D$ is amenable if for
any $\e>0$ and finite dimensional $K$-subspace $1\in E\subset D$, there
exists a finite dimensional $K$-subspace $V\subset D$ such that
$$\dim_K EV < (1+\e) \dim_K V\,.$$
\noindent
All commutative fields and skew fields of finite Gelfand-Kirillov dimension
are amenable. Also, if for a torsion-free amenable group $G$, the group
algebra $KG$ is a domain, then its classical field of fraction is
an amenable skew field.
\begin{propo} \label{amesofic}
All amenable skew fields are sofic.
\end{propo}
\proof It is enough to prove that
for any finite dimensional subspace $1\in F\subset D$ and $\e>0$,
there exists $n\geq 1$
and a $K$-linear map $\tau:D\to\Mat_{n\times n}(K)$ such that
\begin{itemize}
\item $\tau(1)=\Id$
\item $\rk_n(\tau(d))>1-\e$ if $0\neq d\in F$
\item $\rk_n(\tau(fg)-\tau(f)\tau(g))<\e$ for any pair $f,g\in F$.
\end{itemize}
\noindent
First we need a technical lemma.
\begin{lemma}\label{technical}
If $D$ is an amenable skew field, then for any $\delta>0$ and 
finite dimensional subspace $1\in E\subset D$, there exists a pair
of linear subspaces $V_1\subset V\subset D$ such that
$EV_1\subset V$ and
$$\dim_K(V_1)\geq(1-\delta) \dim_K(V)\,.$$ \end{lemma}
\proof
Let $\e>0$ such that $\frac{1}{1+\e}\geq 1-\delta\,.$ By the definition
of amenability, there exists a finite dimensional $K$-subspace
$V_1\subset D$ such that $\dim_K(EV_1)<(1+\e)\dim_K(V_1).$ Now set
$V=EV_1$. Then the pair $V_1\subset V$ satisfies the condition of the
lemma. \qed
\vskip 0.2in
\noindent
Now let $H\subset D$ be the linear subspace spanned by $F\cdot F\,.$
Let $V\subset D$ be an $n$-dimensional linear subspace such that
for some linear subspace $V_1\subset V$, $HV_1\subset V$ and 
$\dim_K(V_1)\geq (1-\e) \dim_K(V)$ hold.
Let $W\subset D$ be a linear subspace complementing $V$ and
let $P:D\to V$ be the $K$-linear projection onto $V$ such that
$$P_{\mid V}=\Id\quad\mbox{and}\quad P_{\mid W}=0\,.$$
Now we define $\tau(d)$ as $P\circ M_d$, where $M_d$ is the left-multiplication
by $d$, Clearly, $\tau:D\to \End_K(V)\cong\Mat_{n\times n}(K)$ is a $K$-linear
map satisfying $\tau(1)=\Id$. If $d\in F$, then
$\Ker(\tau(d))\cap V_1=0\,,$ that is
$\rk_n(\tau(d))\geq 1-\e\,.$ Also, if $f,g\in F$, then
$\tau(fg)=\tau(f)\tau(g)$ restricted on $V_1$. Therefore,
$\rk_n(\tau(fg)-\tau(f)\tau(g))<\e\,.$ \qed

\vskip 0.2in
\section{The free skew field is sofic}
By Theorem 1 of \cite{Elek2006}, if $K\subset E$ is a sub-skew field
of $D$ and $E$ is non-amenable, then $D$ is non-amenable as well. If
$K=\C$ then the free skew field on $r$ generators over $K$ is a non-amenable
skew field. We conjecture that this is the case for all base fields.
There are many ways to define the free skew fields, in this paper
we regard these objects as the skew fields of noncommutative rational
functions. As it follows, we use the approach of \cite{KVV} and \cite{KVV2}.
Let $\cK[z_1,z_2,\dots,z_r]$ be the free algebra over $r$ noncommutative
indeterminates, where $\cK$ be a finite field. 
Formal expressions over $\{z_1,z_2,\dots,z_r\}$ can be defined 
inductively the following way.
\begin{itemize}
\item The elements of  $\cK[z_1,z_2,\dots,z_r]$ are formal expressions.
\item If $R_1$ and $R_2$ are formal expressions, so are $R_1R_2$, $R_1+R_2$
and $(R_1)^{-1}.$
\end{itemize}
\noindent
Let $M_1,M_2,\dots, M_r\in\Mat_{n\times n}(\cK)$, for some $n\geq 1$.
We say that the formal expression $R$ can be evaluated on 
$\{M_1,M_2,\dots, M_r\}$, if all the inverses involved in the 
inductive calculation of $R(M_1,M_2,\dots,M_r)$ exist.
Then we say that $\{M_1,M_2,\dots, M_r\}\in \dom R$. We call a
formal expression $R$ a noncommutative rational expression (nre) if
$\dom R\neq \emptyset$. Two nre's $R$ and $S$ are
equivalent if $\dom R\cap \dom S\neq \emptyset$ and if
$\{M_1,M_2,\dots, M_r\}\in \dom R\cap \dom S$ then
$$R(M_1,M_2,\dots,M_r)=S(M_1,M_2,\dots,M_r)\,.$$
The equivalence classes of nre's are called noncommutative rational
functions and they form the free skew field $\cD_r(\cK)$ over $\cK$ on 
$r$ generators (Proposition 2.2 \cite{KVV}).
Let $T=(T^1_n,T^2_n,\dots,T^r_n)$ be
a $r$-tuple of $n\times n$-matrices with entries
$\{T^i_{jk}\}_{1\leq i\leq r,1\leq j,k\leq n}$ that are commuting
indeterminates. That is, each matrix $T^i_n$ can  be viewed as an element
of the ring $\Mat_{n\times n}(P_{(n)})$, where $P_{(n)}$ is the commutative
 polynomial
algebra over $\cK$ with $rn^2$ variables. The evaluation
$$p\to p(T^1_n,T^2_n,\dots,T^r_n)$$
defines a homomorphism $\rho_n:\cK[z_1,z_2,\dots,z_r]\to \Mat_{n\times n}(P_{(n)})$.
The algebra $\mbox{Im}(\rho_n)=G_n$ is called the algebra of generic matrices.
Let $Q_{(n)}$ be the field of fraction of the algebra $P_{(n)}$. The skew field
$D_n$ is defined as division closure
of $G_n$ in $\Mat_{n\times n}(Q_{(n)})$ 
(Proposition 2.1 \cite{KVV}).
If $R$ is a nre, then $\cN_R$ is defined the
following way. The natural number $n$ is an element of $\cN_R$ if
$$R_n:=R(T_n^1,T_n^2,\dots,T_n^r)$$
can be evaluated inductively as an element of $D_n$.
For any noncommutative rational expression $R$,
there exists $n_R$ such that 
\begin{itemize}
\item If $n\geq n_R$ then $n\in \cN_R$. [Remark 2.3,\cite{KVV2}].
\item If $R$ and $S$ are equivalent rational expressions
 and $n\in \cN_{R}\cap \cN_{S}$, then $R_n=S_n$. 
[Remark 2.6 and Definition 2.8,\cite{KVV2}].
\item For any pair of expressions $R$, $S$ if $n\in \cN_{R}\cap \cN_{S}$, then
$\cN_{RS}=\cN_{R}\cap \cN_S$, $\cN_{R+S}=\cN_{R}\cap \cN_S$
then $(RS)_n=R_nS_{n}, (R+S)_n=R_n+R_m$, whenever $n\in \cN_{R}\cap \cN_S$.
Also, if $n\in \cN_R\cap \cN_{R^{-1}}, R^{-1}_n=(R_n)^{-1}.$ 
[Definition 2.1,\cite{KVV2}].
\end{itemize}
\begin{lemma} \label{utolso1}
If $1\in E\subset \cD_r(\cK)$ is a finite dimensional $\cK$-subspace and
let $F$ be the subspace spanned by $E\cdot E$. Then
there exists $n\geq 1$ and a $\cK$-linear map $\phi:F\to D_n$ such that
$\phi(ab)=\phi(a)\phi(b)$, if $a,b\in E$.
\end{lemma}
\proof
For any element $a\in F$, let us pick a nre $\hat{a}$ that represents $a$.
Choose an integer large enough so that 
$$ n\geq \max\left(\max_{a\in F}\,( n_{\hat{a}})\,\,,
\max_{a,b\in F}\,( n_{\hat{ab}})\,\,,
\max_{a,b\in F}\, (n_{\hat{a+b}}) \right)\,.$$
\noindent
Then, $(\hat{a+b})_n=\hat{a}_n+\hat{b}_n\,.$
Indeed, $\hat{a+b}$ is equivalent to $\hat{a}+\hat{b}$, hence by the basic 
properties above,
$$(\hat{a+b})_n=(\hat{a}+\hat{b})_n=\hat{a}_n+\hat{b}_n\,.$$
\noindent
Similarly,  $(\hat{ab})_n=\hat{a}_n\hat{b}_n\,.$ Therefore, $\phi(a):=\hat{a}_n$
defined a linear map from $F$ into $D_n$ such that $\phi(ab)=\phi(a)\phi(b)$,
whenever $a,b\in E$.\quad \qed
\begin{lemma} \label{utolso2}
$D_n$ is an amenable skew field.
\end{lemma}
\proof First note that $G_n$ is an amenable domain. Indeed, it is a 
subalgebra of
the matrix algebra $\Mat_{n\times n}(P_{(n)})$, hence it has polynomial growth.
 Therefore,
by Proposition 2.2 \cite{Elek2006} its classical field of fractions is an 
amenable skew field.
However, since $D_n$ is the division closure of $G_n$ in a ring, $D_n$ must 
be the
classical field of fractions of $G_n$. \qed
\vskip 0.2in
\noindent
Since amenable skew fields are sofic, by Lemma \ref{utolso1} 
we have the following
proposition.

\begin{propo}
The free skew field on $r$ generators over our finite base field $\cK$
is sofic.
\end{propo}
\proof
It is enough to prove that for any $\cK$-linear subspace
 $1\in E\subset \cD_r(\cK)$ and $\e>0$, there exists $s\geq 1$ and
a unital linear map $\Omega:\cD_r(\cK)\to \Mat_{s\times s}(\cK)$ such that
\begin{equation}\label{E1}
\mbox{For any\,} a,b\in E, \rk_s(\Omega(ab)-\Omega(a)\Omega(b))<\e.
\end{equation}
\begin{equation}\label{E2}
\mbox{For any\,}a\in E, \rk_s(\Omega(a))\geq 1/2\,.
\end{equation}
\noindent
Let $\phi:F\to D_n$ be the map defined in Lemma \ref{utolso1}. Since
$D_n$ is an amenable skew field, there exists $s\geq 1$, and a unital linear map
$\tau:D_n\to \Mat_{s\times s}(\cK)$ such that
\begin{itemize}
\item
For any $c,d\in\phi(F)$, $\rk_s(\tau(cd)-\tau(c)\tau(d))<\e$.
\item
For any $c\in \phi(F)$, $\rk_s(\tau(c))>1/2\,.$
\end{itemize}
\noindent
Now, let $\Omega(a):=\tau\circ\phi(a)$, whenever $a\in F$ (and extend 
$\Omega$ onto $\cD_r(\cK)$ linearly). Clearly, (\ref{E1}) and (\ref{E2}) will
be satisfied. \qed

\section{Limits of linear representations}
Let $\{\theta_k:\F_r\to\Mat_{n_k\times n_k}(\cK)\}^\infty_{k=1}$
be a convergent sequence of linear representations.
Then we can consider the ultraproduct representation
$\theta:\F_r\to \cM_M$ and the extended algebra homomorphism (we denote it with
the same letter) $\theta:\cK\F_r\to \cM_M$. One should notice that
$\cK\F_r/ \Ker(\theta)$ is a sofic algebra.
By definition, for any $n\geq 1$ and $A\in\Mat_{n\times n}(\cK \F_r)$
$$\lim_{k\to\infty} \rk^n_{n_k}(\theta_k(A))=\rk^n_{\cM} (\theta(A))\,.$$
\noindent
Note that $\Mat_{n\times n}(\cM_M)$ is the algebraic ultraproduct of the
matrix algebras $\{\Mat_{n\times n}(\Mat_{n_k\times n_k}(\cK))\}^\infty_{k=1}$ and the
unique extended rank function \cite{Halperin2} on  $\Mat_{n\times n}(\cM_M)$
is exactly the ultralimit of the matrix ranks $\rk^n_{n_k}$.
Now, let us prove Theorem \ref{elso}. It is enough to see that there exists a
von Neumann regular countable subalgebra $S\subset \cM_M$ containing
$\Ima (\theta)$. Indeed, the limit object sought after in the theorem is
the metric closure of $S$ in  $\cM_M$ (recall that the completion of
a rank regular algebra is a continuous algebra \cite{Halperin1}).
If $T$ is a countable subset of  $\cM_M$, then
let $R(T)$ be the $\cK$-algebra generated by $T$. It is easy to see that
$R(T)$ is still countable. Let $X:\cM_M\to \cM_M$ be a function
such that for any $a\in \cM_M$,
$$aX(a)a=a\,.$$
Let $R_1=R(\Ima(\theta)\cup X(\Ima(\theta))$ and inductively,
let
$R_{n+1}=R(R_n\cup X(R_n)$. 
Then the ring $S=\cup^\infty_{n=1}R_n$
is a countable von Neumann regular algebra containing $\Ima(\theta)$.
This finishes the proof of Theorem \ref{elso}. \qed

\vskip 0.2in
\noindent
Now let us suppose that the convergent sequence of representations 
$\{\theta_k\}^\infty_{k=1}$ satisfies the Atiyah condition, that is for
any $n\geq 1$ and $A\in \Mat_{n\times n}(\cK \F_r)$
$$\lim_{k\to\infty} \rk^n_{n_k}(\theta_k(A))\in\Z\,.$$
\begin{propo} \label{linnell}
The division closure $D$ of $\Ima(\theta)$ in the algebra $\cM_M$ is
a skew field. That is, $\theta:\F_r\to D$ is the
limit of the sequence $\{\theta_k\}^\infty_{k=1}$. 
\end{propo}
\noindent
\proof We use an idea of Linnell \cite{Linnell}.
Let $p\neq 0$ be an element of the rational closure of $\Ima (\theta)$ in
$\cM_M$. By Lemma \ref{cramer}, we have matrices $Q,A,B$ such that
all the entries of $Q$ are from $\Ima(\theta)$ and
$$Q=A \begin{pmatrix} 1 & 0 \\ 0 & p \end{pmatrix}\,B\,,$$
where $A$ and $B$ are invertible, and $\rk^n_{\cM}(Q)$ is an integer.
Since the matrix rank on $\Mat_{n\times n}(\cM_M)$ is the ultralimit of
the matrix ranks $\rk^n_{n_k}$ (or by \cite{Halperin2})
$$n-1<\rk^n_{\cM} \begin{pmatrix} 1 & 0 \\ 0 & p \end{pmatrix} \leq n\,.$$
Since 
$$\rk^n_{\cM} \begin{pmatrix} 1 & 0 \\ 0 & p \end{pmatrix} =
\rk^n_{\cM} (Q)$$
by the integrality condition,
$\rk^n_{\cM} \begin{pmatrix} 1 & 0 \\ 0 & p \end{pmatrix} =n$.
Thus, $p$ is invertible. By Lemma \ref{division}, the division
closure of $\Ima(\theta)$ in $\cM_M$ is contained by the
rational closure. Therefore, each nonzero element of the division closure
is invertible. Hence, the division closure $D$ is a skew field. \qed

\vskip 0.2in
\noindent
Now suppose that for some skew field
$D'$, $\theta':\F_r \to D'$ is another limit for
the sequence $\{\theta_k\}^\infty_{k=1}$.
Then for any $n\geq 1$ and $A\in \Mat_{n\times n}(\cK \F_r)$
$$\rk^n_D(\theta(A))=\rk^n_{D'}(\theta(A))\,.$$
Therefore, by Corollary \ref{korolla}, there exists an isomorphism
$\pi: D\to D'$ such that $\pi\circ\theta=\theta'$. This finishes the
proof of Theorem \ref{masodik}. \qed

\vskip 0.2in
\noindent
{\bf Remark.} Note that the proof of Proposition \ref{linnell} shows
that if $\theta':\F_r\to S$ is a limit of the sequence
 $\{\theta_k\}^\infty_{k=1}$, where
$S$ is a continuous algebra, then the
division closure of $\Ima(\theta')$ in $S$ is still isomorphic to $D$.

\section{Linear tilings}
In this section we prove a key technical result of our paper.
Let $D$ be a countable skew field over $\cK$ with $\cK$-basis
$\{1=r_1,r_2,\dots\}$ and let $\phi:D\to\Mat_{n\times n}(\cK)\cong\End_{\cK}(\cK^n)$
be a unit preserving linear map.
We say that $x\in \cK^n$ is $i$-good with respect to $\phi$ if
$$x\in\Ker(\phi(ab)-\phi(a)\phi(b)),$$
whenever $a,b\in \Span\{r_1,r_2,\dots,r_i\}\,.$
The $i$-good elements form the $\cK$-subspace $G^{i,\phi}_n\subset \cK^n\,.$
By definition, if $\{\phi_k:D\to \Mat_{n\times n}(\cK)\}^\infty_{k=1}$
is a sofic representation, then for any fixed $i\geq 1$,
$$\lim_{k\to\infty} \frac{\dim_{\cK} G^{i,\phi}_{n_k}}{n_k}=1\,.$$
\begin{defin}\label{igood}
A unit preserving linear map $\phi:D\to \Mat_{n\times n}(\cK)$ is an $i$-good
map if
$$ \frac{\dim_{\cK} G^{i,\phi}_{n}}{n}\geq 1-\frac{1}{i}\,.$$
\end{defin}

\noindent
Let $\phi:D\to \Mat_{n\times n}(\cK)$ and let $1\in F \subset D$ and 
$H\subset \cK^n$
be finite
dimensional linear subspaces. We call a subset $T\subset \cK^n$ 
a set of $(i,F,H)$-centers with respect to $\phi$ if
\begin{itemize}
\item For any $x\in T$, $\phi(F)(x)\subset H$ is a 
$\dim_{\cK}(F)$-dimensional
$\cK$-subspace.
\item If $x\in T$, then for any $0\neq f\in F$, $\phi(f)(x)$ is $i$-good.
\item The subspaces $\{\phi(F)(x)\}_{x\in T}$ are independent.
\end{itemize}

\noindent
We say that $\phi$ has an $(F,H,i,\delta)$-tiling if there exists a set
of $(i,F,H)$-centers $T$ for $\phi$ such that
$$|T|\dim_{\cK}(F)\geq (1-\delta)n\,.$$

\begin{theorem} \label{tiling}
Let $D$ be a countable skew field over the base field $\cK$, with
basis $\{1=r_1,r_2\dots\}$. Then 
for any finite dimensional subspace $1\in F\subset D$ and $\delta>0$,
we have a positive constant $N_{F,\delta}$ such that if $i,n\geq N_{F,\delta}$,
$\phi:D\to \Mat_{n\times n} (\cK)$ is an $i$-good unit preserving linear
map and $\dim_{\cK}H\geq (1-\frac{1}{i})n,$ then $\phi$
has an $(F,H,i,\delta)$-tiling.
\end{theorem}
\proof Let $\phi:D\to \Mat_{n\times n}(\cK)$ be an $i$-good unit preserving linear
map, the exact values of $i$ and $n$ will be given later.
First note, that if $i$ is larger than some constant
$N^1_F$, then
$(F\backslash\{0\})^{-1},F\subset \Span\{r_1,r_2,\dots,r_i\}\,.$
Note that if
$$0\neq x\in \Ker(\phi(f^{-1})\phi(f)-1)$$
for any $0\neq f\in F$, then
$\dim_{\cK}\phi(F)(x)=\dim_{\cK}(F)\,.$ 
Thus, if $i\geq N^1_F$ and $0\neq x\in G^{i,\phi}_n$,
then $\dim_{\cK} \phi(F)(x)=\dim_{\cK}(F)\,.$
Let $$A_{F,i}=\{x\in \cK^n\,\mid\,\phi(f)(x)\in G^{i,\phi}_n\cap H\,,
\mbox{for any $f\in F$}\}\,.$$

\noindent
It is easy to see that there exists some constant $N^2_{F,\delta}$ such that
if $i\geq N^2_{F,\delta}$ then
\begin{equation} \label{271}
\dim_{\cK}(A_{F,i})\geq (1-\frac{\delta}{4}) n\,.
\end{equation}
\noindent
and for any $0\neq f \in F$,
\begin{equation} \label{272}
\dim_{\cK} \Ker \phi(f)=n-\dim_{\cK} \Ima \phi(f)\leq \frac{\delta}{3} n\,.
\end{equation}
\noindent
Finally, let $N^3_{F,_\delta}>0$ such that
if $n> N^3_{F,_\delta}$, then $|F|\leq |\cK|^{\frac{\delta}{3} n}$ and let
$N_{F,\delta}=\max(N^1_F,N^2_{F,\delta},N^3_{F,_\delta})$.
For $0\neq f\in F$ and $v\in \cK^n$, let
$L(f,v)$ denote the set of points $y$
such that $\phi(f)(y)=v$. Then we have the estimate
\begin{equation} \label{273}
|L(f,v)|\leq |\Ker \phi(f)|\leq |\cK|^{\frac{\delta}{3}n}\,.
\end{equation}
\noindent
Let $T$ be a maximal set of $(i,F,H)$-centers for $\phi$.
We need to prove that
$$|T|\dim_{\cK}(F)\geq (1-\delta) n\,.$$
Let $V$ be the span of the subspace
$\cup_{x\in T}\phi(F)(x)$.
Then,
$$|V|=|\cK|^{\dim_{\cK}(V)}=|\cK|^{|T|\dim_{\cK}(F)}\,.$$
Assume that $|V|<|\cK|^{(1-\delta)n}$.
Now, suppose that $i,n\geq N_{F,\delta}$. Then
$$| \bigcup_{v\in V} \bigcup_{f\in F\backslash\{0\}} L(f,v)|<
|V| |F| |\cK|^{\frac{\delta}{3}n}\leq |\cK|^{(1-\frac{\delta}{3})n}\leq
|A_{F,i}|\,.$$
\noindent
Therefore, there exists $x\in A_{F,i}$ such that
$\phi(F)(x)\cap V=0\,.$ Hence, $\{x\}\cup T$ is a set of
$(i,F,H)$-centers for $\phi$, leading to a contradiction. \qed
\section{Convergent sequences and sofic approximations}
\label{convergent}
The goal of this section is to prove
the following theorem.
\begin{theorem} \label{negyedik}
Let $D$ be a countable skew field over $\cK$ and $\phi:\F_r\to D$ be a 
generating homomorphism. Then $\phi$ is the limit of a convergent
sequence of finite dimensional representations satisfying the Atiyah condition
if and only if $D$ is sofic.
\end{theorem}
\begin{propo}\label{p611}
Let $\phi:\F_r\to D$ as above, where $D$ is sofic. Then $\phi$ is the limit
of a convergent sequence of finite dimensional representations.
\end{propo}
\proof  Let $\{\psi_k:D\to\Mat_{n_k\times n_k}(\cK)\}^\infty_{k=1}$ be
a sofic approximation sequence for the skew field $D$.
\begin{lemma}
For any $B\in\Mat_{n\times n}(D)$
$$\lim_{k\to\infty} \rk^n_{n_k}(\psi_k(B))=\rk^n_D(B)\,,$$
where $\rk^n_{n_k}$ is the matrix rank on $\Mat_{n\times n}(\mnk)\,.$
\end{lemma}
\proof
By taking a subsequence, we may suppose that \\
$\lim_{k\to\infty} \rk^n_{n_k}(\psi_k(B))$ exists for all matrices $B$.
Then, $$\lim_{k\to\infty} \rk^n_{n_k}(\psi_k(B))=\rk^n_{\cM}(\psi(B))\,,$$
where $\psi:D\to M_{\cM}$ is the ultraproduct embedding.
Thus, $\hat{\rk}(B):=\rk^n_{\cM}(\psi(B))$ defines
a rank function on $\mnd$. Since there exists only one rank function
on matrix rings,
$$\lim_{k\to\infty} \rk^n_{n_k}(\psi_k(B))=\rk^n_D(B)\,.\quad\qed$$

\vskip 0.2in
Let $\hp_k:=\psi_k\circ\phi$ be a linear map. By the lemma above,
for any $n\geq 1$ and matrix $A\in\Mat_{n\times n} (\cK\F_r)\,.$
$$\lim_{k\to\infty} \rk^n_{n_k}(\hp_k(A))=\rk^n_D(\phi(A))\,.$$
Note however, that $\hp_k$ does not necessarily 
define a linear representation of
$\F_r$. However, we have the following lemma.
\begin{lemma} \label{l611}
Let $\{\hp_k\}^\infty_{k=1}$ be the maps as above. \\
Let $\{\phi_k:\F_r\to\mnk\}^\infty_{k=1}$ be linear
representations such that
for any generator $\gamma_i$ of the free group
$$\lim_{k\to\infty} \rk_{n_k}(\phi_k(\gamma_i)-\hp_k(\gamma_i))=0\,.$$
Then for any $n\geq 1$ and $A\in\Mat_{n\times n}(\cK \F_r)$
$$\lim_{k\to\infty} \rk^n_{n_k}(\phi_k(A)-\hp_k(A))=0\,.$$
\end{lemma}
\proof
Clearly, it is enough to show that for any $a\in\cK\F_r$
\begin{equation}
\label{e611}
\lim_{k\to\infty} \rk_{n_k}(\phi_k(a)-\hp_k(a))=0
\end{equation}
First we prove (\ref{e611}) in a special case.
\begin{lemma} \label{inverse}
$$\lim_{k\to\infty} \rk_{n_k}(\phi_k(\gamma_i^{-1})-\hp_k(\gamma_i^{-1}))=0\,.$$
\end{lemma}
\proof
By soficity,
$$\lim_{k\to\infty} \rk_{n_k}(\hp_k(\gamma_i^{-1})\hp_k(\gamma_i)-1)=0\,.$$
\noindent
Hence by our assumption,
$$\lim_{k\to\infty} \rk_{n_k}(\hp_k(\gamma_i^{-1})\phi_k(\gamma_i)-1)=0\,.$$
\noindent
Thus
$$\lim_{k\to\infty} \rk_{n_k}\left((\hp_k(\gamma_i^{-1})
-\phi_k(\gamma_i^{-1}))\phi_k(\gamma_i)\right)=0\,.$$
\noindent
Since $\phi_k(\gamma_i)$ is an invertible element for all $k$,
the lemma follows. \qed
\vskip 0.2in
\noindent
Now suppose that for some $w_1,w_2\in\F_r$
$$\lim_{k\to\infty} \rk_{n_k}(\phi_k(w_1)-\hp_k(w_1))=0\quad
\mbox{and}\quad \lim_{k\to\infty} \rk_{n_k}(\phi_k(w_2)-\hp_k(w_2))=0\,.$$
\noindent
By soficity,
$$\lim_{k\to\infty} \rk_{n_k}(\hp_k(w_1w_2)-\hp_k(w_1)\hp_k(w_2))=0\,.$$
Since
$$\phi_k(a)\phi_k(b)-\hp_k(a)\hp_k(b)=(\phi_k(a)-\hp_k(a))\phi_k(b)-\hp_k(a)
(\hp_k(b)-\phi_k(b))$$
we have that
$$\lim_{k\to\infty} \rk_{n_k}(\phi_k(w_1w_2)-\hp_k(w_1w_2))=0\,.\quad$$
\noindent
Therefore by induction, for any $w\in\F_r$
$$\lim_{k\to\infty} \rk_{n_k}(\phi_k(w)-\hp_k(w))=0\,.$$
Now (\ref{e611}) follows easily. \qed
\vskip 0.2in
\noindent
We finish the proof of Proposition \ref{p611}\,. Observe that \\
$\lim_{k\to\infty} \rk_{n_k} \hp_k(\gamma_i)=1$ for all the generators, hence
we have invertible elements
$a^k_i\in \mnk$ such that
$$\lim_{k\to\infty}\rk_{n_k}(\hp_k(\gamma_i)-a^k_i)=0\,.$$
Now, let us define $\phi_k:\F_r\to\mnk$ by
setting $\phi_k(\gamma_i)=a^k_i\,.$ Then the proposition immediately follows
from Lemma \ref{l611}\,.\qed
\vskip 0.2in
\begin{propo} \label{p64}
Let $\{\theta_k:\F_r\to\Mat_{n_k\times n_k}(\cK)\}^\infty_{k=1}$ be
a convergent sequence of linear representations satisfying the Atiyah condition.
Suppose that for some skew field $D$ and generating map
$\theta:\F_r\to D$, $\theta$ is the limit of $\{\theta_k\}^\infty_{k=1}$.
Then $D$ is sofic.
\end{propo}
We will prove a stronger version of Proposition \ref{p64} that will
be used in the next section.
We call $a\in k\F_r$ an element of length at most $l$, if all
the non-vanishing terms of $a=\sum k_i w_i$ have length (as reduced words)
at most $l$. A matrix $A\in\Mat_{s\times s}(\cK \F_r)$ is
of length at most $l$, if all the entries of $A$ have length at most
 $l$. Now, let $\theta:\F_r \to D$ be a generating homomorphism,
where $D$ is a skew field and $\{1=r_1,r_2,\dots\}$ is a $\cK$-basis for $D$.
Let $\rho:\F_r\to \Mat_{s\times s}(\cK)$ be a linear representation.
We say that a linear map $\phi: D\to \Mat_{s\times s}(\cK)$ is an
$(m,\delta)$-approximate extension of $\rho$ if
\begin{itemize}
\item $\phi$ is an $m$-good  map (see Definition \ref{igood}).
\item For any element $a\in \cK \F_r$ of length at most $m$
$$\rk_s(\phi(\theta(a))-\rho(a))<\delta\,.$$
\end{itemize}
\begin{propo} \label{appro}
Let $\theta:\F_r\to D$ be a linear representation into a countable skew field
$D$ such that $\Ima(\theta)$ generates $D$. Let $\{1=r_1, r_2, \dots\}$
be a $\cK$-basis of $D$. Then for any $m\geq 1$ and $\delta>0$  there exists
a constant $l_{m,\delta}$ such that if for a linear representation
$\rho:\F_r\to \Mat_{s\times s}(\cK)$,
\begin{equation}
\left | \rk^n_s(\rho(A))-\rk^n_D(\theta(A))\right |<\frac{1}{l_{m,\delta}}
\end{equation}
whenever $A\in\Mat_{n\times n}(\cK \F_r), n\leq l_{m,\delta}$ is a 
matrix of length at most $l_{m,\delta}$, then
there exists a $\cK$-linear unit preserving map $\phi:D\to \Mat_{s\times s}(\cK)$
that is an $(m,\delta)$-approximate extension of $\rho$.
\end{propo}
\proof Suppose that the Proposition does not hold for some pair $m,\delta$.
Then there exists a sequence
 $\{\theta_k:\F_r\to \mnk\}^\infty_{k=1}$ of linear representations converging
to $\theta$ such that none of the $\theta_k$'s have 
$(m,\delta)$-approximate extension onto
$D$. Consider that ultraproduct map $\hh:\F_r\to \cM_M$. 
By Proposition \ref{linnell} and Theorem \ref{masodik},
we can extend $\hh$ onto an embedding $\phi:D\to \cM_M$ (that is
$\phi\circ \theta=\hh$). From now on, we follow the proof and the notation
of Proposition \ref{soficity}. For $d\in D$, let
$$\phi(d)=[\{\phi_k(d)\}^\infty_{j=1}]\,,$$ where
$\{\phi_k:D\to \Mat_{n_k\times n_k}(\cK)\}^\infty_{k=1}$ are unital linear
maps.
As observed in the
Proposition \ref{soficity}, 
$$\{k\in\N\mid\, \mbox{$\phi_k$ is $m$-good}\}\in\omega$$
Observe that for any $a\in\cK \F_r$
$$\lim_{\omega}\rk_{n_k}(\phi_k(\theta(a))-\theta_k(a))=0\,.$$
\noindent
Therefore,
$$\{k\in\N\mid\,\rk_{n_k}(\phi_k(\theta(a))-\theta_k(a))<\delta, $$
$$ \mbox{for any $a\in\cK \F_r$ of length at most $m$}\}\in\omega.$$
\noindent
Hence,
$$\{k\in\N\mid\,\phi_k \,\mbox{is a $(m,\delta)$-extension}\}\in\omega$$
\noindent
leading to a contradiction. \qed
\vskip 0.2in
\noindent
{\bf Remark.} One should note that for a domain $R$ (provided it is not
an Ore domain) it is possible to have many non-isomorphic skew fields
with epic embeddings $\phi:R\to D$.
In fact, according to our knowledge, 
there is no finitely generated skew field $D$ countable dimensional
over $\cK$ for which epic embeddings $\phi:\cK \F_r\to D$ known not to exist.
In \cite{Herbera}, infinitely many examples of different epic
embeddings of $\theta:\cK\F_r\to Q$ are given, where the skew fields $Q$ are
the quotient fields of  
certain amenable domains. 
Since all these skew fields $Q$ are amenable, they are sofic, hence by
our result above these $\theta$'s are limits of finite dimensional 
representations.
We cannot make the difference between these embeddings using only the
ranks of group algebra elements. This observation
shows why the use of matrices in 
Definition \ref{crucial} is crucial.

\section{Amenable limit fields}
The goal of this section is to prove the first part of Theorem \ref{harmadik}.
\begin{propo}
Let $\{\theta_k:\F_r \to \mnk\}^\infty_{k=1}$, $n_k\to \infty$ be 
a convergent sequence of representations satisfying the Atiyah condition.
Let $\theta:\F_r\to D$ be a limit representation of $\{\theta_k\}^\infty_{k=1}$,
where $D$ is an amenable skew field and $\Ima(\theta)$ generates $D$.
Then $\{\theta_k\}^\infty_{k=1}$ is a hyperfinite sequence.
\end{propo}
\proof
Let $\{1=r_1,r_2,\dots\}$ be a $\cK$-basis for $D$. 
By Proposition \ref{appro}, we have a sofic approximation sequence
$\{\phi_k:D\to\mnk\}^\infty_{k=1}$ such that for any $a\in\cK \F_r$
$$\lim_{k\to\infty} \rk_{n_k}(\phi_k(\theta(a))-\theta_k(a))=0\,.$$

\noindent
Let $\epsilon>0$ and choose $\delta>0$ in such a way that
$(1-\delta)^2>1-\e$ and $(1-\delta)^{-1}\leq 1+\e$. 
By Lemma \ref{technical}, we have finite
dimensional $\cK$-subspaces $F_1\subset F\subset D$ such that
$\theta(\gamma_s)F_1\subset F$ holds for any generator $\gamma_s$ and
$$\dim_{\cK}(F_1)>(1-\delta) \dim_{\cK} (F)\,.$$
Now let $N_{F,\delta}>0$ be the constant in Theorem \ref{tiling}.
Let $i\geq N_{F,\delta}$ be an integer such that
$$\cup_{s=1}^r \theta(\gamma_s)\cup F\subset \Span\{r_1,r_2,\dots,r_i\}\,.$$
\noindent
By definition, there exists $q\geq 1$ such that if $k\geq q$, then
\begin{itemize}
\item $\phi_k$ is $i$-good.
\item $\dim_{\cK} H_k>(1-\frac{1}{i})n_k\,,$ where
$$H_k=\cap^r_{s=1}\{x\in \cK^{n_k}\,\mid\,\phi_k(\theta(\gamma_s))(x)=
\theta_k(\gamma_s)(x)\}\,.$$
\end{itemize}
\noindent
By Theorem \ref{tiling}, if $k>q$, then
$\phi_k$ has an $(F,H_k,i,\delta)$-tiling. Let $T_k$ be the set of
centers of the tiling above.
For $x\in T_k$, let
$$V_x=\{\phi_k(F_1)(x)\}\,.$$
By our assumptions, 
\begin{itemize}
\item For any $x\in T_k$, $\dim_{\cK}(V_x+\sum^r_{s=1} \theta_k(\gamma_s) V_x)
\leq \dim_{\cK}(F)\,.$
\item $\sum_{x\in T_k} \dim_{\cK}(V_x)\geq (1-\delta)^2 n_k$\,.
\item The subspaces $\{W_x=V_x+\sum^r_{s=1} \theta_k(\gamma_s) V_x\}_{x\in T_k}$
are independent.
\end{itemize}
\noindent
Hence, $\{\theta_k\}^\infty_{k=1}$ is a hyperfinite sequence. Indeed,
$K_\e$ can be chosen as $\dim_{\cK}(F)\,.$ \qed

\section{Non-amenable limit fields}
The goal of this section is to finish the
proof of Theorem \ref{harmadik}, by proving 
the following proposition.
\begin{propo} \label{feltetel}
Let $\{\theta_k:\F_r\to\mnk\}^\infty_{k=1}$
be a convergent sequence of finite dimensional
representations satisfying the Atiyah condition. Suppose that
the generating map
$\theta:\F_r\to D$ is a limit of $\{\theta_k\}^\infty_{k=1}$, where
$D$ is a non-amenable skew field. Then $\{\theta_k\}^\infty_{k=1}$ 
is not hyperfinite.
\end{propo}
\noindent
The proof of the proposition will be given in several steps.
\noindent
Let $\theta:\F_r\to D$ be a generating map, where $D$ is a $\cK$-skew field
with basis \\$\{1=r_1,r_2,\dots\}$.
\noindent
Let $\{\rho_k:D\to \mnk\cong\End_{\cK}(\cK^{n_k})\}^\infty_{k=1}$ be
a sequence of unit preserving linear maps such that some non-trivial
$\cK$-subspaces $L_k\subset \cK^{n_k}$ are fixed with uniform bound
$$\dim_{\cK} L_k\leq C\in\N\,\quad \mbox{ for any $k\geq 1$.}$$
\begin{propo} \label{p11.2}
Suppose that $D,\theta,\{\rho_k\}^\infty_{k=1}$ are as above, $\delta>0$,
and for any $k\geq 1$
\begin{itemize}
\item $\dim_{\cK}(L_k+ \sum^r_{i=1} \rho_k(\theta(\gamma_i))L_k)\leq (1+\delta)
\dim_{\cK}L_k\,.$
\item $\rho_k(ab)(x)=\rho_k(a)\rho_k(b)(x)$, if
$x\in L_k$ and \\ $a,b\in\Span\{r_1,r_2,\dots, r_k\}\,.$
\end{itemize}
Then there exists some integer $m\geq 1$ and a finite dimensional
$\cK$-linear subspace $L\subset D^m$ such that
$$\dim_{\cK}(L+\sum^r_{i=1} \theta(\gamma_i)L)\leq (1+\delta)\dim_{\cK}(L)\,.$$
\end{propo}
\proof Again, let $\omega$ be a nonprincipal ultrafilter on the natural numbers.
If $\{V_k\}^\infty_{k=1}$ are finite dimensional $\cK$-linear vectorspaces, then
their ultraproduct $V=\prod_\omega V_k$ is defined the following way. 
Let $Z\subset \prod^\infty_{k=1} V_k$ be the subspace of sequences
 $\{x_k\}^\infty_{k=1}$ such that
$$\{k\,\mid\,x_k=0\}\in\omega \,.$$
\noindent
Then $V=\prod_\omega V_k:=\prod^\infty_{k=1} V_k/Z$. Observe, that if
$\{W_k\subset V_k\}^\infty_{k=1}$ is a sequence of subspaces, then
$\prod_\omega W_k\subset\prod_\omega V_k\,.$ Also, if 
$\{\zeta_k:D\to \End_{\cK}(V_k)\}^\infty_{k=1}$ is a sequence of 
linear maps, the ultraproduct
map $\zeta:D\to \End_{\cK}(V)$ is defined as 
$\zeta(d)(x)=[\{\zeta_k(x_k)\}^\infty_{k=1}]$, where $x=[\{x_k\}^\infty_{k=1}]\in V$.
\begin{lemma} \label{l11.2}
$L=\prod_\omega L_k$ is a non-trivial finite dimensional subspace of
$K=\prod_\omega \cK^{n_k}\,.$ Also, $\dim_{\cK}(L)=t$, where
$\{k\,\mid\, \dim_{\cK}(L_k)=t\}\in\omega\,.$
\end{lemma}
\proof
Let $\{a^k_1, a^k_2,\dots,a^k_C\}$ be a $\cK$-generator system for $L_k$.
Let $x=[\{x_k\}^\infty_{k=1}]\in L\,.$ Then by finiteness, there
exist elements $\{\lambda_i\in\cK\}^C_{i=1}$ such that
$$\{k\mid\,\sum^C_{i=1} \lambda_i a^k_i=x_k\}\in\omega\,.$$
Therefore, $x=\sum_{i=1}^C \lambda_i a_i$, where $a_i=[\{a_i^k\}^\infty_{k=1}\}]\in
L\,.$\quad\qed

\begin{lemma} \label{l11.3}
The ultraproduct of the finite dimensional spaces
$\{L_k+\sum^r_{i=1}\rho_k(\theta(\gamma_i))(L_k)\}$ is
$L+\sum^r_{i=1}\rho(\theta(\gamma_i))(L)$.
\end{lemma}
\proof All the elements of $L+\sum^r_{i=1}\rho(\theta(\gamma_i))(L)$
can be written in the form of
$$x_0+\sum^r_{i=1}\rho(\theta(\gamma_i))(x_i)\,,$$
where $\{x_0,x_1,x_2,\dots,x_r\}\subset V\,.$ Hence
$$\prod_\omega \{L_k+\sum^r_{i=1}\rho_k(\theta(\gamma_i))(L_k)\}\supset
L+\sum^r_{i=1}\rho(\theta(\gamma_i))(L)\,.$$
On the other hand, all the elements of
$\prod_\omega \{L_k+\sum^r_{i=1}\rho_k(\theta(\gamma_i))(L_k)\}$ can be written
as $$[\{x^k_0+\sum^r_{i=1}\rho_k(\theta)(\gamma_i))x_i^k\}^\infty_{k=1}]=
[\{x^k_0\}^\infty_{k=1}]+\sum^r_{i=1}
\rho(\theta)(\gamma_i))[\{x_i^k\}^\infty_{k=1}]\,.$$
Therefore
$$\prod_\omega \{L_k+\sum^r_{i=1}\rho_k(\theta(\gamma_i))(L_k)\}\subset
L+\sum^r_{i=1}\rho(\theta(\gamma_i))(L)\,.$$
\vskip 0.1in
\noindent
By our conditions, if $k$ is large enough, then
\begin{itemize}
\item $\rho_k(b)\rho_k(c)(x)=\rho_k(bc)(x)$
\item $\rho_k(a)\rho_k(bc)(x)=\rho_k(abc)(x)$
\item $\rho_k(ab)\rho_k(c)(x)=\rho_k(abc)(x)$
\end{itemize}
whenever $x\in L_k$ and $a,b,c\in D$.
Hence for the ultraproduct map $\rho$,
\begin{equation} \label{e616}
\rho(ab)\rho(c)(x)=\rho(a)\rho(b)\rho(c)(x)\,,
\end{equation}
whenever $x\in L$. 

\noindent
Now, we finish the proof of Proposition \ref{p11.2}.
Define the $\cK$-vectorspace $T$, by 
$$T:=\rho(D)(L)\subset K\,.$$
Then by (\ref{e616}), we have an embedding
$$\psi:D\to\End_{\cK}(T)\,$$
defined by $\psi(d)(z)=\sum^t_{i=1}\rho(dd_i)l_i\,,$ where
$\{l_1,l_2,\dots,l_t\}$ is a $\cK$-basis of $L$ and 
$z=\sum^t_{i=1}\rho(d_i)l_i\,.$ Hence, $T$ is a left $D$-vectorspace,
with generating system $\{l_1,l_2,\dots,l_i\}$. Also, 
by Lemma \ref{l11.2} and Lemma \ref{l11.3}, $L$ is finite dimensional and 
for $L\subset T\cong D^m$,
$$\dim_{\cK}(L+\sum^r_{i=1} \theta(\gamma_i)L)\leq (1+\delta)\dim_{\cK}(L)\,.$$
\noindent
This finishes the proof of Proposition \ref{p11.2}.\qed
\vskip 0.2in
\noindent
Recall \cite{Elek2006}(Proposition 3.1), that if $D$ is a countable non-amenable
skew field over $\cK$, then
there exist elements $d_1,d_2,\dots, d_l$ and $\e>0$ such
that for any $m\geq 1$ and finite dimensional $\cK$-subspace $W\subset D^m$,
\begin{equation} \label{vege}\frac{\dim_{\cK} (W+\sum^l_{i=1} d_iW)}{\dim_{\cK} W}>1+\e\,.
\end{equation}
Now let $\theta:\F_r\to D$ be a generating map.
\begin{lemma} \label{group}
Let $\theta:\F_r\to D$ be a generating map, where $D$ is non-amenable. 
Then there exists $\delta>0$ such
that for any $m\geq 1$ and finite dimensional $\cK$-subspace $W\subset D^m$,
$$\frac{\dim_{\cK} (W+\sum^r_{i=1}\theta(\gamma_i)W)}{\dim_{\cK} W}>1+\delta\,,$$
where $\{\gamma_i\}^r_{i=1}$ is the standard generator system for $\F_r$.
\end{lemma}
\proof Suppose that such $\delta>0$ does not exist. Then, we have
a sequence of finite dimensional $\cK$-subspaces
$\{W_j\subset D^{n_j}\}^\infty_{j=1}$ such that
$$\lim_{j\to\infty}
\frac{\dim_{\cK} (W_j+\sum^r_{i=1}\theta(\gamma_i)W_j)}{\dim_{\cK} W_j}=1\,.$$
\noindent
Let $S\subset D$ be the set of elements $s$ in $D$ such that
$$\lim_{j\to\infty}
\frac{\dim_{\cK} (W_j+ sW_j)}{\dim_{\cK} W_j}=1\,.$$
\begin{lemma}
$S$ is the division closure of $\Ima(\cK\F_r)$, that is $S=D$.
\end{lemma}
\proof
Clearly, if $a,b\in S$, then $a+b\in S$ and $ab\in S$.
We need to show that if $0\neq a\in S$ then $a^{-1}\in S$.
Let $a\in S$. Since $a(a^{-1}W_j+W_j)=(W_j+aW_j)$, we get that
$$\lim_{j\to\infty}
\frac{\dim_{\cK} (W_j+ a^{-1}W_j)}{\dim_{\cK} W_j}=1\,,$$
therefore $a^{-1}\in S$.\qed
\vskip 0.2in
\noindent
Thus by (\ref{vege}), Lemma \ref{group} follows. \qed
\vskip 0.2in
\noindent
Now we finish the proof of Proposition \ref{feltetel}. Let \\
$\{\theta_k:\F_r\to \mnk\}^\infty_{k=1}$ be a convergent, hyperfinite sequence
of finite dimensional representations satisfying the Atiyah condition.
Let $\theta:\F_r\to D$ be the limit map, where $D$ is a non-amenable
skew field and $\theta$ is a generating map. Let $\delta>0$ be the constant in 
Lemma \ref{group}.
Let $V^k_1, V^k_2,\dots, V^k_{t_k}\subset \cK^{n_k}$ be $\cK$-linear
subspaces such that
\begin{itemize}
\item For any $1\leq j\leq t_k$, $\dck(V^k_j)\leq K_\delta$.
\item $\{V^k_j+\sum^r_{i=1}\theta_k(\gamma_i)V^k_j\}^{t_k}_{j=1}$
 are independent subspaces such that
$$\dim_{\cK}(V^k_j+\sum^r_{i=1}\theta_k(\gamma_i)V^k_j)<(1+\delta)
\dim_{\cK}(V^k_j)\,.$$
\item $\sum^{t_k}_{j=1}\dck(V^k_j)\geq (1-\delta) n_k\,.$
\end{itemize}
\vskip 0.2in
\noindent

We say that the above subspaces $V^k_m$ and $V^l_n$ are equivalent if
there exists a linear isomorphism
$$\zeta:(V^k_m+\sum^r_{i=1}\theta_k(\gamma_i)V^k_m)\to
(V^l_n+\sum^r_{i=1}\theta_k(\gamma_i)V^l_n)$$
such that $\zeta(V^k_m)=V^l_n$ and
$$\zeta(\theta_k(\gamma_i))(x)=\theta_j(\gamma_i)(\zeta(x))$$
for any generator $\gamma_i$ and $x\in V^k_m\,.$
By the finiteness of the base field and the uniform dimension bound, there
are only finitely many equivalence classes. Hence, by taking a subsequence
we can assume that there exists a constant $\tau>0$ such that for each
$k\geq 1$ there are elements $V^k_1,V^k_2,\dots,V^k_{s_k}$ 
of the above subspaces
with the following properties:
\begin{itemize}
\item All the $V^k_j$'s are equivalent.
\item For each $k\geq 1$, $\sum^{s_k}_{j=1}\dim_{\cK}(V^k_j)\geq \tau n_k\,.$
\end{itemize}
\noindent
For $1\leq j \leq s_k$, let $\zeta^k_j:(V^k_1+\sum^r_{i=1}\theta_k(\gamma_i)V^k_1)\to
(V^k_j+\sum^r_{i=1}\theta_k(\gamma_i)V^k_j)$ be the isomorphism 
showing the equivalence. 
If $0\neq \ul=\{\lambda_1,\lambda_2,\dots,\lambda_{s_k}\}
\subset \cK^{s_k}$
is a $s_k$-tuple of elements of $\cK$, then
$\zeta_{\ul}:\sum^{s_k}_{j=1} \lambda_j \zeta^k_j$
defines an isomorphism from $(V^k_1+\sum^r_{i=1}\theta_k(\gamma_i)V^k_1)$ to
$(W^k_{\ul}+\sum^r_{i=1}\theta_k(\gamma_i)W^k_{\ul})$, where
\begin{itemize}
\item $W^k_{\ul}\subset \cK^{n_k}$ is a $\cK$-subspace of
 dimension $\dim_{\cK}(V^k_1)$
\item $\zeta_{\ul}(\theta_k(\gamma_i)(x))=\theta_k(\gamma_i)\zeta_{\ul}(x)$, 
for all $x\in V^k_1$.
\end{itemize}
\begin{lemma}
Let $\{H_k\subset \cK^{n_k}\}^\infty_{k=1}$ be a sequence of subspaces
such that
$\lim_{k\to\infty} \frac{\dim_{\cK} (H_k)}{n_k}=1\,.$
Then if $k$ is large enough, there exists $\ul$ such that
$W^k_{\ul}\subset H_k$ (Note that we cannot assume that for large enough $k$,
$V^k_j\subset H_k$ for some $j\geq 1$).
\end{lemma}
\proof
Let $y^k_1,y^k_2,\dots,y_t^k,t\leq K_\e$ be a $\cK$-basis for $V^k_1$. We define the
linear map $\iota^k_j:\cK^{s_k}\to \cK^{n_k}$ by
$$\iota^k_j(\ul)=\zeta_{\ul}(y^k_j)\,.$$
Since $\{V^k_1,V^k_2,\dots, V^k_{s_k}\}$ are independent subspaces, $\iota^k_j$ is always
an embedding.
Let $$M^k_j=\{\ul\in \cK^{s_k}\,\mid\,\iota^k_j(\ul)\in H_k\}\,.$$
By our assumption, for any $j\geq 1$,
$$\lim_{k\to\infty} \frac{\dim_{\cK} (M^k_j)}{s_k}=1\,.$$
Hence, if $k$ is large enough, then $\cap^t_{j=1} M^k_j\neq \emptyset$. Therefore,
there exists $\ul\in\cK^{s_k}$ such that
$\zeta_{\ul}(V^k_1)\subset H_k$. \qed
\vskip 0.2in
\noindent
By Proposition \ref{appro}, there exists a sequence of maps \\
$\{\rho_k:D\to\mnk\}^\infty_{k=1}$ and subspaces 
$\{H_k\subset \cK^{n_k}\}^\infty_{k=1}$ such that
\begin{itemize}
\item For any $j\geq 1$, there exists $k_j$ such that
$$\rho_k(ab)(x)=\rho_k(a)\rho_k(b)(x)\,,$$
if $k\geq k_j$, $a,b\in\Span\{r_1,r_2,\dots,r_j\}$ and $x\in H_k$.
\item $\rho_k(\theta(\gamma_i))(x)=\theta_k(\gamma_i)(x)\,,$ if
$x\in H_k$ and $\gamma_i$ is a generator of $\F_r$.
\item $\lim_{k\to\infty}\frac{\dim_{\cK} (H_k)}{n_k}=1\,.$
\end{itemize}
\noindent
Hence, the sequence of maps $\{\rho_k\}^\infty_{k=1}$ and subspaces $L_k=W^k_{\ul}$
satisfy the condition of Proposition \ref{p11.2}. Therefore,
there exists $m\geq 1$ and a finite $\cK$-dimensional subspace 
$L\subset D^m$ such that
$$\frac{\dim_{\cK} (L+\sum^r_{i=1}\theta(\gamma_i)L)}{\dim_{\cK} L}\leq 1+\delta\,,$$
in contradiction with the statement of Lemma \ref{group}. 
This finishes the proof of
Proposition \ref{feltetel}. \qed

\vskip 0.3in
\noindent
Gabor Elek \\
Lancaster University \\
g.elek@lancaster.ac.uk
\end{document}